\documentclass[a4paper,12pt,twoside,reqno]{amsart}

\usepackage[utf8]{inputenc}
\usepackage[plainpages=false,pdfpagelabels=true]{hyperref}
\usepackage{amssymb,amsthm}
\usepackage[margin=1in]{geometry}

\newtheorem{Satz}{Theorem}[section]
\newtheorem{Prop}[Satz]{Proposition}

\newtheorem{Thm}[Satz]{Theorem}

\theoremstyle{definition}

\newtheorem{Bem}[Satz]{Remark}

\newcommand{\tr}{\operatorname{Tr}}

\newcommand{\dv}{\text{ }dV}
\newcommand{\s}{\mathbb{S}}
\newcommand{\abs}[1]{\vert #1\vert}
\newcommand{\N}{\mathbb{N}}
\newcommand{\Z}{\mathbb{Z}}
\newcommand{\R}{\mathbb{R}}
\newcommand{\ag}{\mathfrak{g}}
\newcommand{\ah}{\mathfrak{h}}
\newcommand{\nor}{\mathrm{nor}}
\newcommand{\SO}{\mathrm{SO}}

\newcommand{\SU}{\mathrm{SU}}
\allowdisplaybreaks[1]

\renewcommand{\epsilon}{\varepsilon}

\numberwithin{equation}{section}

\title[On the equivariant stability of harmonic self-maps]{On the equivariant stability of harmonic self-maps of cohomogeneity one manifolds}
\author{Volker Branding}
\address{University of Vienna, Faculty of Mathematics\\
Oskar-Morgenstern-Platz 1, 1090 Vienna, Austria\\}
\email{volker.branding@univie.ac.at}

\author{Anna Siffert}
\address{WWU M\"unster, Mathematisches Institut\\
Einsteinstr. 62\\
48149 M\" unster\\
Germany}
\email{anna.siffert@wwu.de}

\date{\today}
\thanks{The first author gratefully acknowledges the support of the Austrian Science Fund (FWF) through the START-project Y963-N35 of Michael Eichmair and 
the project P30749-N35 “Geometric variational problems from string theory”.
The second author gratefully acknowledges the
supports of
the Deutsche Forschungsgemeinschaft (DFG, German Research Foundation) - Project-ID 427320536 - SFB 1442, as well as Germany's Excellence Strategy EXC 2044 390685587, Mathematics M\"unster: Dynamics-Geometry-Structure.}


\begin{document}

\begin{abstract}
The systematic study of harmonic self-maps on cohomogeneity one manifolds
has recently been initiated by P\"uttmann and the second named author in \cite{MR4000241}.
In this article we investigate the corresponding Jacobi equation describing the equivariant stability of such harmonic self-maps. 
Besides several general statements concerning their equivariant stability we explicitly solve
the Jacobi equation for some harmonic self-maps in the cases of spheres, special orthogonal groups and $\SU(3)$.
In particular, we show by an explicit calculation that for specific cohomogeneity one actions
on the sphere the identity map is equivariantly stable.
\end{abstract} 

\maketitle

\section{Introduction and results}
Harmonic maps are solutions to one of the most studied geometric variational problems and have many 
applications in geometry, analysis and mathematical physics.
They are defined as critical points of the energy of smooth maps \(\phi\) between two Riemannian manifolds
\((M^m,g)\) and \((N^n,h)\) which is defined by
\begin{align}
\label{energy}
E(\phi)=\int_M|d\phi|^2\dv.
\end{align}
Its critical points are characterized by the vanishing of the so-called \emph{tension field}
\begin{align}
\label{tension}
0=\tau(\phi):=\tr_g\nabla d\phi,\qquad \tau(\phi)\in\Gamma(\phi^\ast TN),
\end{align}
which are precisely \emph{harmonic maps}.

The harmonic map equation \eqref{tension} is a second order semilinear elliptic partial differential equation.
Due to its non-linear nature
it is a challenging mathematical endeavor to prove the existence of non-trivial solutions.
In the case of both \(M,N\) being closed and \(N\) having non-positive curvature Eells and Sampson
were able to establish their famous existence theorem which ensures 
that every homotopy class of maps contains a harmonic map
under the aforementioned assumptions \cite{MR164306}.

If the target manifold has positive curvature 
it is substantially more difficult to prove the existence of harmonic maps.
To approach this problem 
it is natural
to study the existence of harmonic maps in the case that both manifolds have a sufficient amount of symmetry,
as described in the books
\cite{MR716320,MR1242555}.
Aiming in this direction harmonic maps between cohomogeneity one manifolds were first studied by Urakawa \cite{MR1214054}.
This work has later been extended by P\"uttmann and the second author in \cite{MR4000241}.
The problem of finding a harmonic map between cohomogeneity one manifolds
reduces to solving a singular boundary value problem for an ordinary differential equation of second order.
This ordinary differential equation simplifies further
if one assumes that both domain and target manifold are the same
leading to the notion of \emph{harmonic self-maps}.

The study of harmonic self-maps of spheres with the round metric was initiated by Bizo\'{n} and Chmaj \cite{MR1369057,MR1436833},
see also the subsequent work of Corlette and Wald \cite{MR1810946}.
Extending the results for harmonic map between cohomogeneity one manifolds 
the second author constructed families of harmonic self-maps of spheres \cite{MR3427685}
and of \(\mbox{SU}(3)\) in \cite{MR3745872}.

\medskip

One important property that characterizes the qualitative behavior of a given harmonic map is its stability.
To understand the latter one has to calculate the second variation of the energy \eqref{energy} and evaluate
the resulting expression at a critical point. If a given harmonic map is stable, then there does not exist
a second harmonic map ``nearby'', meaning that the critical points of \eqref{energy} are isolated.

Soon after Eells and Sampson established their existence result for harmonic maps assuming that the
target manifold has negative curvature, Hartman showed that the condition of negative curvature
ensures that such harmonic maps are stable \cite{MR214004}.

On the other hand, in the case of a spherical domain equipped with the round metric and \(m\geq 3\),
Xin showed that harmonic maps are unstable \cite{MR587168}. Soon thereafter
Leung \cite{MR673586} proved that harmonic maps to spheres are unstable in general as well.
In the same spirit of ideas Ohnita was able to extend these results to the case of harmonic maps
from and to symmetric spaces \cite{MR843811}.

\smallskip

In this article we study the \emph{equivariant stability} of harmonic self-maps on cohomogeneity one
manifolds. Here, equivariant stability refers to the fact that we only consider variations
which are invariant under the cohomogeneity one action. We will often see the phenomenon that
a harmonic map will be unstable in general but can be equivariantly stable as we are only
allowing for a special class of directions in the formula for the second variation of the energy \eqref{energy}.

Similar results have already been obtained in the case of biharmonic maps which represent a fourth order generalization of harmonic maps, see \cite{MR4265170} for an introduction to this topic.
The equivariant stability of rotationally biharmonic maps between models was investigated in \cite[Section 5]{MR3357596}.
The normal stability of biharmonic hypersurfaces in spheres has been studied in \cite{ou2021stability}.

In this article we derive the Jacobi operator associated to harmonic self-maps of compact cohomogeneity one manifolds describing their equivariant stability. Throughout this manuscript we assume that the orbit space of the cohomogeneity one manifold is isometric to $[0,L]$, see Section\,\ref{prelim} for more details.
Using the Jacobi operator's explicit form, we prove that its eigenvalues $\lambda_j$ satisfy the asymptotic
behavior
\begin{align*}
\lambda_j=\frac{\pi^2}{L^2}j^2+O(j)
\end{align*}
for \(j\to\infty\).

Further, we will explicitly solve
the Jacobi equation for some harmonic self-maps in the cases of spheres, special orthogonal groups and $\SU(3)$.
In particular, we show by an explicit calculation that for specific cohomogeneity one actions
on the sphere and on special orthogonal groups the identity map is equivariantly stable. 
For example, we will prove the following theorem:

\begin{Thm}
\label{stability-thm-cohom1}
Let $(g,m_0,m_1)$ be one of the following pairs
\begin{align*}
&(2,m_0,m_1), (3,1,1), (3,2,2), (3,4,4), (3,8,8), (4,m_0,1), (4,2,2), (6,1,1), (6,2,2).
\end{align*}
Then the identity map of $\s^{\tfrac{g(m_0+m_1)}{2}+1}$ is a stable harmonic map
with respect to equivariant variations.
\end{Thm}

Throughout this article we make use of the summation convention, i.e. we tacitly sum over repeated
indices.

\medskip

\noindent{\textbf{Organization}:}
This manuscript is structured as follows.
In Section\,\ref{prelim} we provide preliminaries on harmonic maps between cohomogeneity one manifolds and on the stability of harmonic maps.
The equivariant stability of harmonic maps between cohomogeneity one manifolds is discussed in Section\,\ref{stability-cohom1}, where we in particular provide Theorem\,\ref{stability-thm-cohom1}.
In the last section, Section\,\ref{explicit}, we study the equivariant stability of some explicitly given equivariant harmonic self-maps of cohomogeneity one manifolds and explicitly determine the spectra of the corresponding Jacobi operators.

\section{Preliminaries}
\label{prelim}
In Subsection\,\ref{prelim-harmonic} we give a brief introduction to harmonic maps between cohomogeneity one manifolds.
Afterwards, in Subsection\,\ref{prelim-stability}, we recall various facts concerning the 
stability of harmonic maps which are related to the results of this article.

\subsection{Harmonic maps between cohomogeneity one manifolds}
\label{prelim-harmonic}

In this subsection we introduce relevant notation and results from the theory of equivariant harmonic self-maps of compact cohomogeneity one manifolds. The main source is \cite{MR4000241}.

\smallskip

In what follows let $M$ be a Riemannian manifold endowed with an isometric action $G\times M\to M$ of a compact Lie group $G$.
We further assume that the orbit space $M/G$ is isometric to a closed interval $[0,L]$ and that the Weyl group $W$ of the action is finite, i.e. we are considering specific cohomogeneity one actions on $M$. 
In this setting P\"uttmann and the second named author studied when
equivariant {\em $(k,r)$-maps}, i.e., maps of the form
\begin{align}
\label{def-kr-map}
  g\cdot\gamma(t)\mapsto g\cdot\gamma(r(t))
\end{align}
where $r: [0,L] \to \R$ is a smooth function with $r(0) = 0$ and $r(L) = kL$, are harmonic \cite{MR4000241}.
Here, $\gamma$ denotes a fixed unit speed normal geodesic with $\gamma(0)$ being contained in one of the non-principal orbits. 
Further, the integer $k$ is of the form $j\abs{W}/2+1$ with $j\in 2\Z$. 
We want to mention that for specific group actions also odd integers $j$ might be allowed.
The Brouwer degree of a $(k,r)$-map is given by $k$ if the codimensions of both non-principal orbits are odd, and by $1$ otherwise,
see \cite{MR2480860} for more details.

In the following, let $Q$ be a fixed biinvariant metric on $G$ and denote the orthonormal complement of the Lie algebra $\ah$ of the principal isotropy group $H$ in $\ag$ by $\mathfrak{n}$.
Then, the metric endomorphism $P_t : \mathfrak{n} \to \mathfrak{n}$ is defined by
\begin{gather*}
  Q(X, P_t\cdot Y) = \langle X^{\ast},Y^{\ast} \rangle_{\vert\gamma(t)},
\end{gather*}
where $Z^{\ast}$ denotes the action field of $Z\in\mathfrak{n}$.
In addition, we define \(B_+: \mathfrak{n} \times \mathfrak{n}\to \mathfrak{g}\)
as follows
\begin{align*}
B_+^{r(t)}(X,Y)=\frac{1}{2}([X,P_{r(t)}\cdot Y]-[P_{r(t)}\cdot X,Y]).
\end{align*}

The tension field of a \((k,r)\)-map as defined in \eqref{def-kr-map}
is given by the following expression
\begin{align*}
\tau_{\vert\gamma(t)} = 
\underbrace{\ddot r(t) + \tfrac{1}{2}\dot r(t) \tr P_t^{-1}\dot P_t
   - \tfrac{1}{2} \tr P_t^{-1} \dot P_{r(t)}}_{:=\tau^{\nor}_{\vert\gamma(t)} }
+\underbrace{\sum_{j=1}^{m-1}P^{-1}_{r(t)}B_+^{r(t)}(e_j,e_j)}_{:=\tau^{\tan}_{\vert\gamma(t)}}. 
\end{align*}

Hence, a \((k,r)\)-map is harmonic if and only if
\begin{align}
\label{tension-field-kr-map}
\tau_{\vert\gamma(t)}=\tau^{\nor}_{\vert\gamma(t)}+\tau^{\tan}_{\vert\gamma(t)}=0.
\end{align}
We realize that both the tangent and the normal part of the tension field have to vanish.

For the specific actions we deal with in this paper, the tangential component of the tension field vanishes such that the construction of harmonic maps is reduced to the construction of solutions to the ordinary differential equation
\begin{equation}
\label{ode-harmonic}
\ddot r(t) + \tfrac{1}{2}\dot r(t) \tr P_t^{-1}\dot P_t
    - \tfrac{1}{2} \tr P_t^{-1} \dot P_{r(t)}=0.
\end{equation}

\smallskip

Below we give further details on the particular cases of cohomogeneity one actions on spheres, on special orthogonal groups and on $\SU(3)$.

\subsubsection{Cohomogeneity one actions on spheres}
\label{cohom-sphere}
By proving that each cohomogeneity one action on a sphere is orbit equivalent to the isotropy representation of a Riemannian symmetric space of rank $2$, Hsiang and Lawson \cite{MR298593} classified such actions.
The orbits of any isometric cohomogeneity one action $G\times \s^{n+1} \to \s^{n+1}$ of a compact Lie group $G$
on the sphere $\s^{n+1}$ yield an isoparametric foliation of the sphere.
 The latter is a family of parallel hypersurfaces with constant principal curvatures together with two focal submanifolds which together foliate the sphere.
 Takagi and Takahashi \cite{MR0334094}
  determined the number $g$ of possible distinct principal curvatures of the orbits (i.e. the isoparametric hypersurfaces) and their multiplicities $m_0,\ldots,m_{g-1}$. They thus proved
\begin{gather*}
  m_0 = m_2 = \ldots = m_{g-2} \quad \text{and} \quad m_1 = m_3 = \ldots = m_{g-1}.
\end{gather*}
Therefore, we have $n = \frac{m_0+m_1}{2} g$. 
Later, M\"unzner \cite{MR583825} showed that this holds true for all isoparametric foliations of spheres, i.e. also for those not stemming from group actions.
Up to ordering of $m_0$ and $m_1$ there exist only actions with the following $(g,m_0,m_1)$:
\begin{multline*}
 (1,m,m), (2,m_0,m_1), (3,1,1), (3,2,2), (3,4,4), (3,8,8),\\
   (4,m_0,1), (4,2,2), (4,2,2\ell+1), (4,4,4\ell+3), (4,4,5), (4,6,9), (6,1,1), (6,2,2),
\end{multline*}
where \(\ell\in\N_+\).

Harmonic self-maps of $\s^{n+1}$
can be characterized as critical points of 
\begin{align*}
E(r)=c\int_0^\frac{\pi}{g}\big(\dot r(t)^2+\sum_{i=0}^{g-1}m_i\tfrac{\sin^2(r(t)-i\tfrac{\pi}{g})}{\sin^2(t-i\tfrac{\pi}{g})}
\big)\Pi_{i=0}^{g-1}\sin^{m_i}(t-i\tfrac{\pi}{g})dt,
\end{align*} 
where \(c\) denotes a positive constant.
The critical points are those which are solutions to the ordinary differential equation
\begin{align}
\label{euler-lagrange-isoparametric-sphere}
  0 =
  4\sin^{2}(gt) \cdot \ddot r(t) + \bigl( g(m_0+m_1)\sin(2gt) + 2g(m_0-m_1)\sin(gt) \bigr)\dot r(t)\\\notag
  - g(g -2)\sin(2(r(t)-t)) \bigl( m_0+m_1 + (m_0-m_1)\cos(gt) \bigr)\\\notag
  -2g \sin\bigl(2(r(t)-t)+gt\bigr) \bigl( (m_0+m_1)\cos(gt) +m_0-m_1\bigr)
\end{align}
for functions $r: \,]0,\frac{\pi}{g}[\, \to \R$ with
\begin{gather}
\label{bv-sphere}
  \lim_{t\rightarrow 0}r(t)=0 \quad \text{and}\quad \lim_{t\rightarrow \tfrac{\pi}{g}}r(t)=k\tfrac{\pi}{g}, 
\end{gather}
where $k\in\Z$.
The above boundary value problem \eqref{euler-lagrange-isoparametric-sphere}, \eqref{bv-sphere}
will be referred to as \emph{$(g,m_0,m_1,k)$-boundary value problem}.
When the multiplicities $m_0$ and $m_1$ coincide, we will refer to this problem as
 \emph{$(g,m,k)$-boundary value problem}.
It can be easily seen that \(r(t)=t\) is an explicit solution of the $(g,m_0,m_1,1)$-boundary value problem.
It is clear that this solution corresponds to the identity map
of the Riemannian manifold \(M\).
Further, the function \(r(t)=(1-g)t\) is an 
explicit solution of the  $(g,m,1-g)$-boundary value problem.
In the following we will refer to these particular solutions as \emph{linear solutions}.
It has been shown in \cite{MR4000241} that the $(g,m_0,m_1,k)$-boundary value problem admits linear solutions for $k=1$ and $k=1-g$, $m_0=m_1$ only and that the above mentioned solutions exhaust all linear solutions.

\smallskip

In \cite{MR1436833} Bizo\'{n} and Chmaj proved that each of the $(1,m,0)$-boundary value problems and the $(1,m,1)$-boundary value problems admits infinitely many solutions
 provided that $2\leq m\leq5$. 
For each choice of $m$ with $2\leq m\leq5$, we subsume the solutions with $k=0$ and $k=1$ in a 
countable family $r_n$, $n\in\N$, of solutions which are labeled by a nodal number $n$, namely by the number of intersections with $\tfrac{\pi}{2}$.
Bizo\'{n} and Chmaj also proved
that these solutions (more precisely, the solutions different from the identity map or its negative) do no longer exist for \(m\geq 6\).

\begin{Bem}
Note that Bizo\'{n} and Chmaj \cite{MR1436833} use a different notation than we do: the role of $m$ is played by what they call $k$, 
which satisfies
 $k=m+1$. Further, instead of $r(t)$ they deal with the shifted function $h(t)=r(t)-\tfrac{\pi}{2}$.
\end{Bem}

In \cite{MR3427685} the second named author provided an analogous result for $g=2$, namely
she proved that each of the $(2,m_0,m_1,1)$-boundary value problems admits infinitely many solutions
assuming that $2\leq m_0\leq5$.  
For each choice of $m_0$ with $2\leq m_0\leq5$, we thus obtain a 
countable family $r_n$, $n\in\N$, of solutions which are labeled by a nodal number $n$, namely by the number of intersections with $\tfrac{\pi}{2}$.
It was also shown in  \cite{MR3427685}
that these solutions (more precisely, the solutions different from the identity map or its negative) do no longer exist for \(m_0\geq 6\).

\subsubsection{Cohomogeneity one actions on special orthogonal groups}
Any isometric cohomogeneity one action $G\times \s^{n+1}\to \s^{n+1}$ on a sphere can be lifted to an isometric cohomogeneity one action of $G\times\SO(n+1)$ on $\SO(n+2)$. 
If we denote the data of the cohomogeneity one action on the sphere by $(g,m_0,m_1)$, then we have $n=(m_0+m_1)g/2$. 
It was proved in \cite{MR4000241} that a
solution of the $(2g,m_0,m_1,k)$-boundary value problem yields a harmonic self-map of $\SO(n+2)$.
We can conclude that harmonic self-maps of $\SO(n+2)$
are associated with the solutions of the ordinary differential equation

\begin{align}
\label{euler-lagrange-isoparametric-so}
  0 =
  4\sin^{2} (2gt) \cdot \ddot r(t) + \bigl( 2g(m_0+m_1)\sin(4gt) + 4g(m_0-m_1)\sin (2gt) \bigr)\dot r(t)\\\notag
  - 4g(g -1)\sin(2(r(t)-t)) \bigl( m_0+m_1 + (m_0-m_1)\cos(2gt) \bigr)\\\notag
  -4g \sin\bigl(2(r(t)-t )+2gt\bigr) \bigl( (m_0+m_1)\cos(2gt) +m_0-m_1\bigr)
\end{align}
for functions $r: \,]0,\frac{\pi}{2g}[\, \to \R$ with
\begin{gather}
\label{bv}
  \lim_{t\rightarrow 0}r(t)=0 \quad \text{and}\quad \lim_{t\rightarrow \tfrac{\pi}{2g}}r(t)=k\tfrac{\pi}{2g}, 
\end{gather}
where $k\in\Z$.
Clearly, $r(t)=t$ is a solution to this boundary value problem with $k=1$.
Further, if $m=m_0=m_1$, then $r(t)=(1-2g)t$ is a solution to this boundary value problem with 
$k=1-2g$.

\subsubsection{Cohomogeneity one action on $\SU(3)$}
For the cohomogeneity one action
\begin{align*}
\SU(3)\times\SU(3)\rightarrow\SU(3),\hspace{1cm} (A,B)\mapsto ABA^{\mbox{tr}},
\end{align*}
it has been shown in \cite{MR3745872} that
\begin{align*}
P_t=4\,\mbox{diag}\left(1,\cos^2(t),\cos^2(t),\sin^2(t/2),\sin^2(t/2),\cos^2(t/2),\cos^2(t/2)\right),
\end{align*}
assuming that $\SU(3)$ is endowed with the metric $\langle A_1,A_2\rangle=\mbox{tr}(A_1\overline{A_2}^{\mbox{tr}})$
and the normal geodesic $\gamma$ is given by
\begin{align*}
\gamma(t)=\left(\begin{smallmatrix}
\cos t&-\sin t&0\\
 \sin t & \cos t&0\\
0&0&1
\end{smallmatrix} \right).
\end{align*}
Further, it was proven that each solution to 
\begin{align*}
\ddot r(t)=-\csc^2(2t)\left(2\sin(4t)\dot r(t)+4\sin^2(t)\sin2r(t)-8\cos^3(t)\cdot\sin r(t)\right),
\end{align*}
which satisfies
$r(0)=0$ and $r(\tfrac{\pi}{2})=(2\ell+1)\tfrac{\pi}{2}, \ell\in\Z$, yields a harmonic self-map of $\SU(3)$.
More precisely, it was shown in \cite{MR3745872} that there exists an infinite, countable sequence of solutions of the above boundary value problem with either $\ell=0$ or $\ell=1$. Numerical investigations indicate that we have $\ell=0$ for these solutions but this has not been proved.

\subsection{Stability of harmonic maps}
\label{prelim-stability}
In this subsection we recall various facts concerning the 
stability of harmonic maps where we follow \cite[Chapter 5]{MR1252178}.

Thus, let \(\phi\colon M\to N\) be a smooth harmonic map.
The second variation of the energy of a map \eqref{energy} evaluated at a critical point
is given by
\begin{align}
\delta^2 E(\phi)(V,W)=\int_M\langle J_\phi(V),W\rangle\dv,\qquad V,W\in\Gamma(\phi^\ast TN).
\end{align}
Here, \(J_\phi\) denotes the \emph{Jacobi operator} which is defined by
\begin{align}
\label{jacobi-operator}
J_\phi(V):=\Delta^\phi V-R^N(V,d\phi(e_i))d\phi(e_i).
\end{align}
In the above formula \(R^N\) represents the Riemannian curvature tensor of the manifold \(N\).
We use the following sign convention for the \emph{rough Laplacian}
\begin{align*}
\Delta^\phi V:=-(\nabla^\phi_{e_i}\nabla^\phi_{e_i}-\nabla^\phi_{\nabla_{e_i}e_i})V,
\end{align*}
where \(\nabla^\phi\) represents the connection on \(\phi^\ast TN\) and \(\{e_i\},i=1,\ldots,m\)
is an orthonormal basis of \(TM\). Note that the Jacobi operator \(J_\phi\) defined in \eqref{jacobi-operator}
is self-adjoint when \(M\) is compact.

We say that a harmonic map \(\phi\) is \emph{stable} if
\begin{align}
\delta^2 E(\phi)(V,V)>0\qquad \textrm{ for all } V\in\Gamma(\phi^\ast TN).
\end{align}
Moreover, a harmonic map \(\phi\) is called \emph{weakly stable} if
\(\delta^2 E(\phi)(V,V)\geq 0\), otherwise 
it is called \emph{unstable}.

From the general spectral theory for elliptic operators over a compact Riemannian manifold
we know that the eigenvalues \(\lambda\) of \(J_\phi\) satisfy
\begin{align*}
\lambda_1(\phi)\leq\lambda_2(\phi)\leq\ldots\leq\lambda_j(\phi)\to \infty.
\end{align*}

The vector space
\begin{align*}
V_\lambda(\phi):=\{\xi\in\Gamma(\phi^\ast TN)\mid J_\phi\xi=\lambda\xi\}\neq\{0\}
\end{align*}
is called \emph{eigenspace} with eigenvalue \(\lambda\).

Moreover, \(\dim V_\lambda(\phi)\) is called the \emph{multiplicity} of \(\lambda\)
and we know from general elliptic theory that \(\dim V_\lambda(\phi)<\infty\). 

In terms of these spectral data we define 
\begin{align*}
\mathrm{index}(\phi)&:=\sum_{\lambda <0}\dim V_\lambda(\phi),\\
\mathrm{nullity}(\phi)&:=\dim V_0(\phi)=\dim\ker (J_\phi).
\end{align*}

It follows directly that a harmonic map \(\phi\) is stable if and only if
\(\lambda_j(\phi)>0\) for \(j\in\N_+\).

\medskip

In the following we will cite four important results concerning the stability of 
harmonic maps which are closely related to the main results of this article,
the first two have already been mentioned in the introduction.

\begin{Satz}[Xin]
\label{theorem-xin}
For \(m\geq 3\) any stable harmonic map from
\(\s^m\) with the round metric into any Riemannian manifold must be a constant map.
\end{Satz}

In order to prove the above result the author make use
of the fact that the sphere admits conformal vector fields.

\begin{Satz}[Leung]
\label{theorem-leung}
Let \((N^n,h)\) be a complete orientable hypersurface of \(\R^{n+1}\).
Denote the principal curvatures by \(\lambda_1,\ldots,\lambda_n\) and 
set \(\lambda^2:=\max\lambda_i^2,i=1,\ldots,n\).
Moreover, let \(K\) be the function that assigns to each point in \(N\)
the minimum of the sectional curvature at that point.
If \(\lambda^2<(n-1)K\), then any stable harmonic map \(\phi\colon M\to N\)
is constant.
\end{Satz}

In particular, this implies that
for \(n\geq 3\) a stable harmonic map from
any Riemannian manifold into \(\s^n\) with the round metric must be a constant map.

The above results have been extended to Riemannian symmetric spaces in \cite{MR843811,MR850121}.

On the other hand, it was shown by Urakawa that by a suitable deformation of the metric on the domain
a harmonic map to the sphere can also be stable \cite[Proposition 7.4]{MR882704}.
For further results on the stability of harmonic maps we refer to \cite[Chapter 5]{MR1252178}.

Moreover, 
let us recall the following result concerning the stability of the identity map on \(\s^{m+1}\),
which was obtained independently by Mazet \cite[Proposition 8]{MR336767} and Smith \cite[Example 2.12]{MR375386}.

\begin{Prop}[Mazet, Smith]
Consider the identity map of \((\s^{m+1},g)\) with \(m\geq 2\) as a solution to the harmonic map equation.
Then the smallest proper eigenvalue of the Jacobi operator is \(1-m\) with multiplicity \(m+2\).
The other proper eigenvalues are strictly positive.
\end{Prop}

Here, proper means that we exclude the eigenspace spanned by the constant function.
The above result follows from calculating the spectrum of the operator
\begin{align}
\label{operator-sphere-general}
L:=\Delta_{\s^{m+1}}-2m\operatorname{Id}
\end{align}
acting on functions. It is well known that 
\begin{equation*}
\operatorname{spec}\Delta_{\s^{m+1}}=\{\lambda_j=j(j+m)\,\lvert\,\,j\geq 0\}
\end{equation*}
and each eigenvalue has multiplicity \(m+2\).

For the sake of completeness we also want to mention the following
result providing a geometric characterization of the stability of the identity map,
which is also due to Smith \cite{MR375386}.

\begin{Satz}[Smith]
Let \((M,g)\) be a compact Einstein manifold with Einstein constant \(\rho\).
Then the following statements hold:
\begin{enumerate}
 \item The identity map of \(M\) is weakly stable if and only if the first eigenvalue
 of the Laplacian acting on smooth functions, denoted by \(\lambda_1(g)\), satisfies
 \begin{align}
  \lambda_1(g)\geq 2\rho.
 \end{align}
\item The nullity of the identity map is given by
\begin{align}
\operatorname{nullity(Id)}=\dim\operatorname{Iso}(M,g)
+\dim\{f\in C^\infty(M)\mid \Delta f=2\rho f\},
\end{align}
where \(\operatorname{Iso}(M,g)\) is the isometry group of \((M,g)\), i.e.
\begin{align*}
\operatorname{Iso}(M,g):=\{\varphi\mid\varphi^\ast g=g\}.
\end{align*}
\end{enumerate}
\end{Satz}

In the following we will study the equivariant stability of harmonic self-maps 
between cohomogeneity one manifolds. Such maps may be stable with respect to equivariant
variations but unstable with respect to general variations. Hence, our results on the equivariant
stability of some particular harmonic maps do not contradict the general theorems that
have been presented in this subsection.

\section{Equivariant stability of harmonic self-maps between cohomogeneity one manifolds}
\label{stability-cohom1}
In the present section we discuss the equivariant stability of equivariant harmonic self-maps of compact cohomogeneity one manifolds $M$. 
We make use of the notation introduced in Subsection\,\ref{prelim-harmonic}.

In order to investigate the equivariant stability of harmonic \((k,r)\)-maps,
which are characterized as solutions of \eqref{tension-field-kr-map},
we consider a variation of \(r(t)\), denoted by \(r_s(t)\), which satisfies
\begin{align*}
\frac{d}{ds}\big|_{s=0}r_s(t)=\xi(t).
\end{align*}
In the following we will calculate the linearization of the normal and the tangential
part of the tension field separately.
The variation of the normal part is given by
\begin{align}
\label{jacobi-normal}
0&=\frac{d}{ds}\big|_{s=0}\big(\tau^{\nor}_{\vert\gamma(t)}\big)  \\
\nonumber&=\frac{d}{ds}\big|_{s=0}\big(\ddot r_s(t) + \tfrac{1}{2}\dot r_s(t) \tr P_t^{-1}\dot P_t
   - \tfrac{1}{2} \tr P_t^{-1} \dot P_{r_s(t)}\big) \\
\nonumber&=\ddot\xi(t)+\tfrac{1}{2}\,\tr(P_t^{-1}\dot P_t)\dot\xi(t)-\tfrac{1}{2}\tr(P_t^{-1}\ddot P_{r(t)})\xi(t),
\end{align}
whereas for the tangential part we find
\begin{align}
\label{jacobi-tangential}
0&=\frac{d}{ds}\big|_{s=0}\big(\tau^{\tan}_{\vert\gamma(t)}\big)  \\
\nonumber&=\frac{d}{ds}\big|_{s=0}\big(\sum_{j=1}^{m-1}P^{-1}_{r_s(t)}B_+^{r_s(t)}(e_j,e_j)\big) \\
\nonumber&=\sum_{j=1}^{m-1}\frac{d}{dt}\big(P^{-1}_{r(t)}B_+^{r(t)}(e_j,e_j)\big)\xi(t).
\end{align}

Hence, we obtain the following

\begin{Prop}
Let $r(t)\colon [0,L]\to\R$ be a solution of \eqref{tension-field-kr-map} and $\phi$ be the associated harmonic self-map on a cohomogeneity one manifolds.
The corresponding Jacobi equation is given by
\begin{align}
\label{jacobi-general}
\ddot\xi(t)+\tfrac{1}{2}\,\tr(P_t^{-1}\dot P_t)\dot\xi(t)-\tfrac{1}{2}\tr(P_t^{-1}\ddot P_{r(t)})\xi(t)
+\sum_{j=1}^{m-1}\frac{d}{dt}\big(P^{-1}_{r(t)}B_+^{r(t)}(e_j,e_j)\big)\xi(t)
+\lambda\xi(t)=0,
\end{align}
where \(\xi\in C^{\infty }_{0}(\left[0,L\right])\).
\end{Prop}
\begin{proof}
This is a direct consequence of the above calculations.
\end{proof}

Note that we can write the Jacobi equation \eqref{jacobi-general} in the following form
\begin{align}
\label{jacobi-sl-form}
\frac{d}{dt}\big(\sqrt{\det P_t}~\dot\xi(t)\big)
+\big(-\tfrac{1}{2}\tr(P_t^{-1}\ddot P_{r(t)})+\sum_{j=1}^{m-1}\frac{d}{dt}\big(P^{-1}_{r(t)}B_+^{r(t)}(e_j,e_j)\big)
\big)&\sqrt{\det P_t}~\xi(t) \\
\nonumber&+\lambda\sqrt{\det P_t}~\xi(t)=0,
\end{align}
which allows us to apply 
the general theory of Sturm-Liouville for one-dimensional eigenvalue problems.
For more details on Sturm-Liouville theory as it is used in this article
we refer to Appendix \ref{appendix-sturm-liouville}.

We get the following general result:
\begin{Satz}
\label{theorem-ev-general}
Let \(r(t)\colon [0,L]\to\R\) be a solution of the equation characterizing
harmonic self-maps \eqref{tension-field-kr-map}.
Then the following statements for the eigenvalue problem of the corresponding Jacobi operator \(\xi(t)\) hold:
\begin{enumerate}
 \item The eigenvalue problem \eqref{jacobi-general} has infinitely many simple eigenvalues
\begin{align*}
\lambda_0<\lambda_1<\lambda_2<\ldots\qquad \lambda_j\to\infty\textrm{ for } j\to\infty.
\end{align*}

 \item The eigenfunction \(\xi_j(t)\) corresponding to the eigenvalue \(\lambda_j\) has exactly \(j\) zeros
in \((0,L)\). Between two zeros of \(\xi_j(t)\) there is exactly one zero of \(\xi_{j+1}(t)\).

\item For \(j\to\infty\) the eigenvalues have the following asymptotic behavior
\begin{align}
\label{asymptotic-ev-jacobi}
\lambda_j=\frac{\pi^2}{L^2}j^2+O(j).
\end{align}
\end{enumerate}
\end{Satz}

\begin{proof}
Performing the substitutions \(z(t)=p(t)=\sqrt{\det P_t}\) and 
\begin{align*}
q(t)=\big(\sum_{j=1}^{m-1}\frac{d}{dt}\big(P^{-1}_{r(t)}B_+^{r(t)}(e_j,e_j)\big)
-\tfrac{1}{2}\tr(P_t^{-1}\ddot P_{r(t)})\big)\sqrt{\det P_t}
\end{align*}
we conclude that the Jacobi equation \eqref{jacobi-sl-form} can be written in the form of a
Sturm-Liouville eigenvalue problem, see Appendix \ref{appendix-sturm-liouville} 
for the precise details.

Moreover, it can be directly seen that the first three conditions of \((SL)\),
which are given in Appendix \ref{appendix-sturm-liouville} in full detail,
are satisfied. 
Regarding the fourth condition of \((SL)\) we note that
the boundary conditions \eqref{sl-boundary-conditions}
are satisfied by choosing \(\alpha_1=\beta_1=0\)
and \(\alpha_2,\beta_2\neq 0\) as \(p(t):=\sqrt{\det P_t}\)
satisfies \(p(0)=p(\frac{\pi}{2g})=0\) and hence the fourth condition of \((SL)\)
also holds true.

The statement \eqref{asymptotic-ev-jacobi} on the asymptotic behavior of the eigenvalues
is a direct consequence of the Weyl asymptotic \eqref{sl-weyl}.
\end{proof}

\section{Explicit calculation of spectra}
\label{explicit}
In this section we study the equivariant stability of some explicitly given equivariant harmonic self-maps of cohomogeneity one manifolds.
These harmonic maps have been provided in \cite{MR4000241,MR3745872}.

We assume that we only deal with those cohomogeneity one actions for which the tangential component of the tension field vanishes trivially, in other words, for which solutions to \eqref{ode-harmonic} induce harmonic maps between the corresponding cohomogeneity one manifolds.
For these kinds of maps also the tangential contribution to the Jacobi equation \eqref{jacobi-tangential} vanishes.
The cohomogeneity one actions on spheres, special orthogonal groups and on $\SU(3)$ discussed in Subsection\,\ref{prelim-harmonic} satisfy this condition.

\subsection{Equivariant stability of harmonic self-maps of the sphere}
Throughout this section we use that
the Jacobi operator corresponding to the Euler-Lagrange equation 
\eqref{euler-lagrange-isoparametric-sphere} is given by
\begin{align}
\label{jacobi-isoparametric-general}
&\ddot\xi(t)+\big(m_0-m_1+(m_0+m_1)\cos(gt)\big)\tfrac{g}{2\sin(gt)}\dot\xi(t)\\
\notag&-\frac{g}{2\sin^2(gt)}\big[(g-2)\cos(2(r(t)-t))\big(m_0+m_1+(m_0-m_1)\cos(gt)\big)\\\notag&+2\cos(2(r(t)-t)+gt)(m_0-m_1+(m_0+m_1)\cos(gt))\big]\xi(t) \\
\nonumber&+\lambda\xi(t)=0.
\end{align}

\subsubsection{The equivariant stability of the linear solution \texorpdfstring{r(t)=t}{}}
\label{stability-rt}
In this subsection we investigate the equivariant stability of the linear solution \(r(t)=t\) of \eqref{euler-lagrange-isoparametric-sphere}. 
Plugging in $r(t)=t$ 
and performing a change of variables 
\(t=\frac{2}{g}\arctan(e^x)\) equation \eqref{jacobi-isoparametric-general} transforms into the following eigenvalue problem
\begin{align}
\label{jacobi-gg-identity}
\xi''(x)&+\tfrac{1}{2}\big(m_0-m_1+(2-(m_0+m_1))\tanh(x)\big) \xi'(x)\\\notag&
-\tfrac{1}{2}(m_0+m_1-(m_0-m_1)\tanh(x))\xi(x)
+\big(\lambda+\frac{m_0+m_1}{g}\big)\frac{\xi(x)}{\cosh^2(x)}=0.
\end{align}
(The constant $\lambda$ in the previous equation differs from $\lambda$ in \eqref{euler-lagrange-isoparametric-sphere} by a factor $g^2$. By abuse of notation we call both constants $\lambda$.)

In order to solve \eqref{jacobi-gg-identity}, we make the ansatz
\begin{align*}
\xi(x)=\frac{1}{\cosh (x)} f(x)
\end{align*}
and obtain
\begin{align}
\label{ode-gg-f}
f''(x)&-\tfrac{1}{2}(-m_0+m_1+(2+m_0+m_1)\tanh(x)) f'(x)\\
\notag&+\tfrac{1}{2g}\big(g(-2+2\lambda-m_0-m_1)+2(m_0+m_1)\big)\tfrac{f(x)}{\cosh^2(x)}=0.
\end{align}

It can be read of directly that
\begin{align*}
f(x)=1,\qquad \lambda=\tfrac{m_0+m_1}{2}+1-\tfrac{m_0+m_1}{g}
\end{align*}
is a solution of \eqref{ode-gg-f}.

In order to find additional solutions of \eqref{ode-gg-f} 
and to solve the spectral problem \eqref{jacobi-gg-identity}
we perform the transformation \(f(x)=u(\tanh(x))\) which gives the equation
\begin{align*}
(1-\tanh^2(x))u''(\tanh(x))&- \tfrac{1}{2}(-m_0+m_1+(6+m_0+m_1)\tanh(x))u'(\tanh(x))
\\\notag&+\tfrac{1}{2g}(g(-2+2\lambda-m_0-m_1)+2(m_0+m_1))u(\tanh(x))=0.
\end{align*}

The above equation is solved by the Jacobi polynomials,
see Appendix \ref{appendix-polynomial-jacobi} for the precise details. We summarize our calculations as follows:

\begin{Prop}
\label{prop-sphere-m0m1}
The spectral problem \eqref{jacobi-gg-identity} describing the equivariant stability of the identity
map, which we parametrize by \(r(x)=\frac{2}{g}\arctan(e^x)\), is solved by
\begin{align}
\label{spectrum-gg-linear1}
\xi_j(x)=\frac{1}{\cosh(x)}P^{(\tfrac{m_1+1}{2},\tfrac{m_0+1}{2})}_{j-1}\big(\tanh(x)\big),\qquad \lambda_j=-\tfrac{m_0+m_1}{g}+j(j+\tfrac{m_0+m_1}{2}),
\end{align}
where \(j\in\,\mathbb{N}_+\).
\end{Prop}

Note that in Proposition \ref{prop-sphere-m0m1} we have used the convention \(j\geq 1\) and normalized
the eigenvalues such that $\lambda_1=-\tfrac{m_0+m_1}{g}+1+\tfrac{m_0+m_1}{2}$.

As an immediate consequence of Proposition\,\ref{prop-sphere-m0m1} we obtain the following theorem.

\begin{Thm}
\label{thm-stable-sphere}
Let $(g,m_0,m_1)$ be one of the following pairs
\begin{align*}
&(1,1,1),(2,m_0,m_1), (3,1,1), (3,2,2), (3,4,4), (3,8,8), (4,m_0,1), (4,2,2), 
(6,1,1), (6,2,2).
\end{align*}
Then the identity map of $\s^{\tfrac{g(m_0+m_1)}{2}+1}$ is equivariantly stable.
\end{Thm}

\begin{Bem}
\label{rem-g1}
\begin{enumerate}
\item Theorem\,\ref{thm-stable-sphere} states that the identity map is equivariantly stable for all possible values of $(g,m_0,m_1)$
except $(1,m)$, $m\in\N_+$ with $m\geq 2$
and except the triples $(4,2,2\ell+1), (4,4,4\ell+3), (4,4,5), (4,6,9)$ where $\ell,m_0,m_1\in\N_+$ with $m_0\leq m_1$.
We would like to point out once more that this 
does not contradict the general instability Theorem \ref{theorem-xin} as we are only considering a special class of variations.
\item For the case $g=1$, Proposition\,\ref{prop-sphere-m0m1} has been established by Bizo\'{n} and Chmaj in \cite{MR1436833}.
 \item Note that for $g=1$ the spectrum  \eqref{spectrum-gg-linear1}
 is precisely the spectrum of the operator \eqref{operator-sphere-general}
which was calculated by abstract methods instead of a direct calculation.
\item Setting $g=1$ and performing the change of variables \(t=2\arctan(e^x)\), \eqref{euler-lagrange-isoparametric-sphere} transforms into 
\begin{align*}
r''(x)-(m-1)\tanh(x)r'(x)-\frac{1}{2}m\sin 2r(x)=0.
\end{align*} 
By differentiating this identity one finds that \(r'(x)\)
solves the eigenvalue problem \eqref{jacobi-gg-identity} with eigenvalue \(\lambda=1-m\).
There is a geometric reason for this fact: The proof of Theorem \ref{theorem-xin} makes
use of a conformal vector field on the sphere which is obtained by projecting
a parallel vector field from \(\R^{m+2}\) onto \(\s^{m+1}\).
As \(\frac{\partial}{\partial x}\) is the generator of conformal transformations
it is clear that \(r'(x)\) solves \eqref{jacobi-gg-identity} with the corresponding eigenvalue.
However, for \(g\geq 2\) the above statement does no longer hold true.
\item 
If we inspect the complete list of triples  $(g,m_0,m_1)$ presented in Subsection \ref{cohom-sphere}
then we realize that the triples
\begin{align*}
(4,2,2\ell+1), (4,4,4\ell+3), (4,4,5), (4,6,9), 
\end{align*}
where $\ell,m_0,m_1\in\N_+$ with $m_0\leq m_1$, do not appear in the statement of Theorem \ref{thm-stable-sphere}. It is currently unknown if the tangential part of the tension field 
\eqref{tension-field-kr-map} vanishes for these triples and hence we cannot make a statement
on the equivarant stability of the identity map in these cases.
\end{enumerate}
\end{Bem}

\subsubsection{The equivariant stability of the linear solution \texorpdfstring{r(t)=(1-g)t}{}}
\label{stability-rg}
As a next step we study the equivariant stability of the linear solution \(r(t)=(1-g)t\) of \eqref{euler-lagrange-isoparametric-sphere}. 
Note that this solutions exists only for $m=m_0=m_1$. 
Performing the change of variables 
\(t=\frac{2}{g}\arctan(e^x)\) equation \eqref{jacobi-isoparametric-general} acquires the following form
\begin{align}
\label{jacobi-gg-linear2}
\xi''(x)-(m-1)\tanh (x) \xi'(x)-m\tanh^2(x)\xi(x)+\big(\lambda+m-2\frac{m}{g}\big)\frac{\xi(x)}{\cosh^2(x)}=0.
\end{align}

Again, we make the ansatz
\begin{align*}
\xi(x)=\frac{1}{\cosh (x)} f(x)
\end{align*}
leading to 
\begin{align}
\label{ode-gg-linear2-f}
f''(x)-(m+1)\tanh(x) f'(x)+\big(\lambda-2\frac{m}{g}+m-1\big)\frac{f(x)}{\cosh^2(x)}=0.
\end{align}

We conclude that 
\begin{align*}
f(x)=1,\qquad \lambda=2\frac{m}{g}-m+1
\end{align*}
is a solution of \eqref{ode-gg-linear2-f}.

Again, to obtain the additional solutions of \eqref{ode-gg-linear2-f} 
and to solve the spectral problem \eqref{jacobi-gg-identity},
we perform the transformation \(f(x)=u(\tanh(x))\) which results in the equation
\begin{align*}
(1-\tanh^2(x))u''(\tanh(x))-&(3+m)\tanh(x) u'(\tanh(x))\\
&+(\lambda-2\frac{m}{g}+m-1)u(\tanh(x))=0.
\end{align*}

Hence, by the same reasoning as in the previous section we find that 
the above equation is solved by the so-called Gegenbauer polynomials,
see Appendix \ref{appendix-polynomial-gegenbauer} for the precise details.

\begin{Prop}
\label{prop-sphere-2}
The spectral problem \eqref{jacobi-gg-linear2}, describing the equivariant stability of the linear solution
parametrized by \(r(x)=(1-g)\frac{2}{g}\arctan(e^x)\), is solved by
\begin{align}
\label{spectrum-gg-linear2}
\xi_j(x)=\frac{1}{\cosh(x)}C^{(\frac{m+2}{2})}_{j-1}\big(\tanh(x)\big),\qquad \lambda_j=j(j+m)+2\frac{m}{g}-2m,
\end{align}
where \(j\in\mathbb{N}_+\).
\end{Prop}

As an immediate consequence of Proposition\,\ref{prop-sphere-2} we obtain the following theorem.

\begin{Thm}
Let $(g,m)$ be one of the following pairs
\begin{align*}
(1,m), (2,m), (3,1), (3,2),(4,1), (4,2), (6,1),
\end{align*}
where $m\in\N_+$.
Then the harmonic self-map of $\s^{mg+1}$, which is associated with the solution $r(t)=(1-g)t$ of \eqref{euler-lagrange-isoparametric-sphere}, 
is equivariantly stable.
\end{Thm}

\begin{Bem}
\begin{enumerate}
\item
In the case \(g=2\) the eigenvalues \eqref{spectrum-gg-linear1} and \eqref{spectrum-gg-linear2}
coincide as one should expect from the explicit form of the linear solutions \(r(t)\).
\item
We would like to point out that the explicit spectra \eqref{spectrum-gg-linear1}, \eqref{spectrum-gg-linear2}
do not contradict Theorems \ref{theorem-xin} and \ref{theorem-leung} as we are dealing with a special kind of stability, that is equivariant stability.
\end{enumerate}
\end{Bem}

\subsection{Equivariant stability of harmonic self-maps of \texorpdfstring{SO(n)}{}}
Throughout this subsection let
$(g,m_0,m_1)$ be one of the following pairs
\begin{align*}
&(1,m_0,m_0),(2,m_0,m_1), (3,1,1), (3,2,2), (3,4,4), (3,8,8), (4,m_0,1), (4,2,2), \\
&(6,1,1), (6,2,2),
\end{align*}
where $\ell,m_0,m_1\in\N_+$ with $m_0\leq m_1$.

\smallskip

The Jacobi operator corresponding to the Euler-Lagrange equation \eqref{euler-lagrange-isoparametric-so}
is given by
\begin{align}
\label{jacobi-son}
&\ddot\xi(t)+\big(m_0-m_1+(m_0+m_1)\cos(gt)\big)\tfrac{g}{\sin(2gt)}\dot\xi(t)\\
\notag&-\frac{g}{\sin^2(2gt)}\big((2g-2)\cos(2(r(t)-t))\big(m_0+m_1+(m_0-m_1)\cos(2gt)\big)\\\notag&+2\cos(2(r(t)-t)+2gt)(m_0-m_1+(m_0+m_1)\cos(2gt))\big)\xi(t) \\
\nonumber&+\lambda\xi(t)=0.
\end{align}

\subsubsection{The equivariant stability of the linear solution \texorpdfstring{r(t)=t}{}}
In this subsection we investigate the equivariant stability of the linear solution \(r(t)=t\) of \eqref{euler-lagrange-isoparametric-so}. 
Plugging in $r(t)=t$ and performing a change of variables 
\(t=\frac{1}{g}\arctan(e^x)\) equation \eqref{jacobi-son} transforms into the following eigenvalue problem

\begin{align}
\label{jacobi-gg-identity-so}
\xi''(x)&+\tfrac{1}{2}\big(m_0-m_1+(2-(m_0+m_1))\tanh(x)\big) \xi'(x)\\\notag&
-\tfrac{1}{2}(m_0+m_1-(m_0-m_1)\tanh(x))\xi
+\big(\lambda+\frac{m_0+m_1}{2g}\big)\frac{\xi(x)}{\cosh^2(x)}=0.
\end{align}

From the considerations in Subsection\,\ref{stability-rt} we get:

\begin{Prop}
\label{prop-so}
The spectral problem \eqref{jacobi-gg-identity} describing the equivariant stability of the identity
map, which we parametrize by \(r(x)=\frac{1}{g}\arctan(e^x)\), is solved by
\begin{align}
\label{spectrum-son}
\xi_j(x)=\frac{1}{\cosh(x)}P^{(\tfrac{m_1+1}{2},\tfrac{m_0+1}{2})}_{j-1}\big(\tanh(x)\big),\qquad \lambda_j=-\tfrac{m_0+m_1}{2g}+j(j+\tfrac{m_0+m_1}{2}),
\end{align}
where \(j\in\,\mathbb{N}_+\).
\end{Prop}

As an immediate consequence of Proposition\,\ref{prop-so} we obtain the following theorem.

\begin{Thm}
\label{thm-stable}
Let $(g,m_0,m_1)$ be one of the following pairs
\begin{align*}
&(1,m_0,m_0),(2,m_0,m_1), (3,1,1), (3,2,2), (3,4,4), (3,8,8), (4,m_0,1), (4,2,2), \\
&(6,1,1), (6,2,2),
\end{align*}
where $\ell,m_0,m_1\in\N_+$ with $m_0\leq m_1$.
Then the identity map of $\SO(n+2)$ is equivariantly stable, where $n=g(m_0+m_1)$.
\end{Thm}

\subsubsection{The equivariant stability of the linear solution \texorpdfstring{r(t)=(1-2g)t}{}}
As a next step we study the equivariant stability of the linear solution \(r(t)=(1-2g)t\) of \eqref{euler-lagrange-isoparametric-so}. 
Note that this solutions exists only for $m=m_0=m_1$. 
Again, performing the change of variables 
\(t=\frac{1}{g}\arctan(e^x)\) equation \eqref{jacobi-son} acquires the following form
\begin{align}
\label{jacobi-so-2}
\xi''(x)-(m-1)\tanh (x) \xi'(x)-m\tanh^2(x)\xi(x)+\big(\lambda+m-\frac{m}{g}\big)\frac{\xi(x)}{\cosh^2(x)}=0.
\end{align}

By the same considerations as in Subsection\,\ref{stability-rg} we get:

\begin{Prop}
\label{prop-so-2}
The spectral problem \eqref{jacobi-so-2}, describing the equivariant stability of the linear solution
parametrized by \(r(x)=\frac{1-2g}{g}\arctan(e^x)\), is solved by
\begin{align}
\label{spectrum-so-linear2}
\xi_j(x)=\frac{1}{\cosh(x)}C^{(\frac{m+2}{2})}_{j-1}\big(\tanh(x)\big),\qquad \lambda_j=j(j+m)+\frac{m}{g}-2m,
\end{align}
where \(j\in\mathbb{N}_+\).
\end{Prop}

As an immediate consequence of Proposition\,\ref{prop-so-2} we obtain the following theorem.

\begin{Thm}
Let $(g,m)$ be one of the following pairs
\begin{align*}
(1,m), (2,1), (2,2), (3,1), (4,1), (6,1),
\end{align*}
where $m\in\N_+$.
Further, we set $n=2gm$.
Then the harmonic self-map of $\SO(n+2)$ which is associated with the solution $r(t)=(1-2g)t$ of \eqref{euler-lagrange-isoparametric-so}, 
is equivariantly stable.
\end{Thm}

\subsection{The second variation of the identity map of \texorpdfstring{\mbox{SU}(3)}{}}
Harmonic self-maps of \(\mbox{SU}(3)\) have been investigated by the second author in \cite{MR3745872}
where the existence of a countably infinite family of harmonic self-maps of \(\mbox{SU}(3)\) with
non-trivial Brouwer degree was established.

The latter can be characterized as critical points of the energy functional
\begin{align}
\label{energy-su3}
E(r)=\int_{-\infty}^\infty\big(r'^2(x)+(1+\tanh(x))\cos^2 r(x)-\sqrt{2}\cos r(x)(1-\tanh(x))^\frac{3}{2}\big)
\frac{dx}{\cosh (x)}.
\end{align}
More precisely, the critical points of \eqref{energy-su3} are those who satisfy
\begin{align}
\label{ode-su3}
r''(x)-\tanh(x) r'(x)-\frac{1+\tanh(x)}{2}\sin 2r(x)
-\frac{1}{\sqrt{2}}(1-\tanh(x))^\frac{3}{2}\sin r(x)=0.
\end{align}
The second variation of \eqref{energy-su3} evaluated at a critical point 
is given by
\begin{align}
\label{second-variation-su3}
\delta^2&E(r)(\xi,\xi)=
\\
&2\int_{-\infty}^\infty\big(\xi'^2(x)-(1+\tanh(x))\cos 2r(x) \xi^2(x)
+\frac{1}{\sqrt{2}}(1-\tanh(x))^\frac{3}{2}\cos r(x) \xi^2(x)\big)\frac{dx}{\cosh(x)}.
\nonumber
\end{align}
Hence, in order to study the equivariant stability of harmonic self-maps of \(SU(3)\) we have to investigate
the following eigenvalue problem
\begin{align}
\label{jacobi-su3-general}
\xi''(x)-\xi'(x)\tanh(x)+&(1+\tanh(x))\cos 2r(x)\xi(x) \\
\nonumber&-\frac{1}{\sqrt{2}}(1-\tanh(x))^\frac{3}{2}\cos r(x)\xi(x)
+\frac{\lambda }{4\cosh^2(x)}\xi(x)=0.
\end{align}
Note that in our coordinates the volume element is given by \(\frac{1}{2}\frac{dx}{\cosh^3(x)}\)
which leads to the factor of \(\frac{1}{\cosh^2(x)}\).

Again, the only solution of \eqref{ode-su3} which is known in closed form
is the identity map \(r_1(x)=\arctan(e^x)\). 
Applying the identity
\begin{align*}
\cos(\arctan e^x)=\frac{1}{\sqrt{1+e^{2x}}}=\frac{1}{\sqrt{2}}\sqrt{1-\tanh(x)}
\end{align*}
we find that \eqref{jacobi-su3-general} simplifies to 
\begin{align}
\label{jacobi-su3-identity}
\xi''(x)-\tanh (x) \xi'(x)-(\frac{1}{2}+\frac{3}{2}\tanh^2(x))\xi(x)+\frac{\lambda}{4\cosh^2(x)}\xi(x)=0.
\end{align}
Now, we make the ansatz
\begin{align*}
\xi(x)=\cosh^2(x) f(x)
\end{align*}
and find 
\begin{align}
\label{ode-su3-f}
\cosh^2(x)f''(x)-3\cosh(x)\sinh(x)f'(x)+(\frac{3}{2}+\frac{\lambda}{4})f(x)=0.
\end{align}
Hence, we can conclude that \(f(x)=1\) and \(\lambda=-6\) is a solution of \eqref{ode-su3-f}.
To obtain the additional solutions we make the ansatz \(f(x)=u\big(\frac{1}{2}\tanh(x)\big)\) which yields
the equation
\begin{align*}
(1-\tanh^2x)u''\big(\frac{1}{2}\tanh(x)\big)
-10\tanh x u'\big(\frac{1}{2}\tanh(x)\big)+(6+\lambda)u\big(\frac{1}{2}\tanh(x)\big)=0.
\end{align*}
This equation is again solved by the Gegenbauer polynomials presented in Appendix \ref{appendix-polynomial}.

We summarize our calculations as follows:
\begin{Prop}
The spectral problem \eqref{jacobi-su3-identity} characterizing the equivariant stability of the identity 
map of \(SU(3)\) parametrized by \(r_1(x)=\arctan(e^x)\) is solved by
\begin{align}
\label{spectrum-su3}
\xi_j(x)=\cosh^2(x)~C^{(\frac{9}{2})}_{j-1}\big(\frac{1}{2}\tanh(x)\big),\qquad \lambda_j=j(j+7)-14,
\end{align}
where \(j\in\N_+\).
\end{Prop}

\begin{Bem}
Recall that \(\dim \mbox{SU}(3)=8\), if we now compare the spectrum of the Jacobi operator for the identity
on \(\s^{m+1}\) in the case of \(m=7\) and $g=1$ given by \eqref{spectrum-gg-linear1}
with the spectrum of the Jacobi operator on \(\mbox{SU}(3)\) given by \eqref{spectrum-su3} we realize
that the eigenvalues coincide, while the eigensections differ slightly.
\end{Bem}

\begin{Bem}
In all explicit calculations that we have carried out in this section,
characterizing the equivariant stability of the identity map, we have seen that 
the eigenvalues \(\lambda_j\) have a growth rate of \(j^2\).
This is consistent with the statement of Theorem \ref{theorem-ev-general}
which describes the general equivariant stability of harmonic self-maps
on cohomogeneity one manifolds.
\end{Bem}

\appendix
\section{Aspects of Sturm-Liouville theory}
\label{appendix-sturm-liouville}

In this appendix we collect a number of results from Sturm-Liouville theory
for ordinary differential equations which are linear and of second order.
For more details on this subject we refer to \cite[Chapter 5.4]{MR2961944}.

We set \(J:=[a,b]\) and consider \(u\colon J\to\R\).
We define the operator \(\mathcal{L}\) as follows
\begin{align*}
\mathcal{L} u(s):=\frac{d}{ds}\big(p(s)\frac{d}{ds}u(s)\big)+q(s)u(s).
\end{align*}
We are interested in the eigenvalue problem
\begin{align}
\label{sl-ev}
\mathcal{L}u(s)+\lambda z(s)u(s)=0
\end{align}
for which we impose the two boundary conditions
\begin{align}
\label{sl-boundary-conditions}
\alpha_1u(a)+\alpha_2p(a)u'(a)=0,\qquad
\beta_1u(b)+\beta_2p(b)u'(b)=0
\end{align}
with \(\alpha_i,\beta_i\in\R,i=1,2\).
Note that \eqref{sl-boundary-conditions} can be read as a linear combination
of Dirichlet and Neumann boundary data for \(u\). 

We say that the conditions \((SL)\) are satisfied if the following conditions hold
\begin{enumerate}
 \item \(p(s)\in C^1(J)\),
 \item \(q(s),z(s)\in C^0(J)\),
 \item \(p(s)>0, z(s)>0\) in \(J\),
 \item \(\alpha_1^2+\alpha_2^2>0\) and \(\beta_1^2+\beta_2^2>0\).
\end{enumerate}

We are now ready to state the following powerful result from Sturm-Liouville theory
which is the basis in the proof of Theorem \ref{theorem-ev-general}.
\begin{Satz}
Consider the eigenvalue problem \eqref{sl-ev}. If the conditions \((SL)\) hold, then
the eigenvalue problem \eqref{sl-ev} has infinitely many simple eigenvalues
\begin{align*}
\lambda_0<\lambda_1<\lambda_2<\ldots\qquad \lambda_j\to\infty\textrm{ for } j\to\infty.
\end{align*}
The eigenfunction \(u_j(s)\) corresponding to the eigenvalue \(\lambda_j\) has exactly \(j\) zeros
in \(]a,b[\). Between two zeros of \(u_j(s)\) there is exactly one zero of \(u_{j+1}(s)\).

Moreover, we have the following asymptotic behavior of the eigenvalues 
(Weyl asymptotics)
for large values of \(j\)
\begin{align}
\label{sl-weyl}
\lambda_j=\pi^2\big(\int_a^b\sqrt{\frac{z(s)}{p(s)}}ds\big)^{-2}j^2+O(j).
\end{align}
\end{Satz}

\section{Aspects of orthogonal polynomials}
\label{appendix-polynomial}
Here, we provide some facts on the specific orthogonal polynomials
which are solutions to the linear second order ordinary differential equations
which appear in the study of the equivariant stability of harmonic self-maps.
For more details on this subject we refer to \cite[Chapter 22]{MR0167642}
and the website \cite[Chapter 18]{NIST:DLMF}.

Consider a second order ordinary differential equation of the from
\begin{align}
\label{ode-general}
A(x)f''(x)+B(x)f'(x)+C(x)f(x)+\lambda_jf(x)=0.
\end{align}

\subsection{Jacobi polynomials}
\label{appendix-polynomial-jacobi}
If 
\begin{align*}
A(x)=1-x^2,\qquad B(x)=\beta-\alpha-(\alpha+\beta+2)x,\qquad C(x)=0, \qquad & \lambda_j=j(j+1+\alpha+\beta)
\end{align*}
then \eqref{ode-general} is solved by the \emph{Jacobi polynomials} \(P_j^{(\alpha,\beta)}(x)\), where \(j\geq 0\)
and \(\alpha,\beta\geq -1\).
The polynomials 
\(P_j^{(\alpha,\beta)}(x)\) are sometimes also called \emph{hypergeometric polynomials} in the literature.

These polynomials can be characterized as
\begin{align*}
 P_j^{(\alpha,\beta)}(x)=\tfrac{(\alpha+1)_n}{n!}{}_2F_1(-n,n+1+\alpha+\beta,\alpha+1,\tfrac{1}{2}(1-x)),
\end{align*}
where ${}_2F_1$ is the hypergeometric function and $(\alpha+1)_n$ denotes the Pochhammer symbol for the rising factorial.

\subsection{Gegenbauer polynomials}
\label{appendix-polynomial-gegenbauer}
If 
\begin{align*}
A(x)=1-x^2,\qquad B(x)=-(2\delta+1)x,\qquad C(x)=0,\qquad \lambda_j=j(j+2\delta)
\end{align*}
then \eqref{ode-general} is solved by the \emph{Gegenbauer polynomials} \(C_j^{(\delta)}(x)\), where \(j\geq 0\)
and \(\delta>-\frac{1}{2},\delta\neq 0\).
The Gegenbauer polynomials are a special case of the Jacobi polynomials introduced above.
The polynomials \(C_j^{(\delta)}(x)\) are sometimes also called \emph{ultraspherical polynomials} in the literature.

Note that we can obtain the Gegenbauer polynomials from the Jacobi polynomials by choosing
\(\alpha=\beta=\delta-\frac{1}{2}\).

For \(0\leq j\leq 3\) these are explicitly given by
\begin{align*}
C_0^{(\delta)}(x)&=1, \\
C_1^{(\delta)}(x)&=2\delta x, \\
C_2^{(\delta)}(x)&=-\delta+2\delta(1+\delta)x^2, \\
C_3^{(\delta)}(x)&=-2\delta(1+\delta)x+\frac{4}{3}\delta(1+\delta)(2+\delta)x^3.
\end{align*}

One possibility of obtaining the higher order Gegenbauer polynomials is the recursion formula
\begin{align*}
C_j^{(\delta)}(x)=\frac{1}{j}\big(2x(j+\delta-1)C_{j-1}^{(\delta)}(x)-(j+2\delta-2)C_{j-2}^{(\delta)}(x)\big),
\qquad j\geq 2.
\end{align*}

\bibliographystyle{plain}
\bibliography{mybib}

\end{document}